# On the Bautin bifurcation for systems of delay differential equations

*Anca–Veronica Ion*
*University of Piteşti,*
*Târgu din Vale, nr.1,*
*Piteşti, Romania*

*Abstract.* For systems of delay differential equations the Hopf bifurcation was investigated by several authors. The problem we solve here is that of the possibility of emergence of a codimension two bifurcation, namely the Bautin bifurcation, for some such systems.

*Keywords:* bifurcation theory, Bautin, delay differential equation

## 1. Introduction

The existence of periodic solutions for evolution equations is of certain interest for both pure and applied mathematicians. Even for bidimensional systems of differential equations the detection of limit cycles by theoretical means is difficult. The bifurcation theory offers a strong tool for finding limit cycles, namely the theory concerning the Hopf bifurcation (when there is a varying parameter)[2], [6]. Several authors studied the Hopf bifurcation for delay differential equations (e.g. [4], [7], [1], [5] ). We are interested to find sufficient conditions for the Bautin bifurcation for a class of such systems.

## 2. Setting of the problem, theoretical framework

Consider a system of the form

$$\dot{x}(t) = A(\alpha)x(t) + B(\alpha)x(t-r) + f(x(t), x(t-r), \alpha), \qquad (1)$$
$$x(s) = \phi(s), \quad s \in [-r, 0], \qquad (2)$$

where $x = (x_1, ..., x_n) \in \mathbf{R}^n$, $\alpha = (\alpha_1, \alpha_2) \in \mathbf{R}^2$, $A(\alpha), B(\alpha)$ are nxn matrices over $\mathbf{R}$, $f = (f_1, ..., f_n)$ is continuously differentiable on its domain of definition, $D \subset \mathbf{R}^{2n+1}$. Moreover, $f(0,0,\alpha) = 0$ and the differential of $f$ in the first two vectorial variables, calculated at $(0,0,\alpha)$ is equal to zero. $\phi$ is an element of the Banach space $B = C([-r,0], \mathbf{R}^n)$ (column vectors). In order to write eq. (1) as a differential equation in a Banach space, the space

$$B_0 = \{\psi : [-r, 0] \to \mathbf{R}^n, \psi \text{ is continuous on } [-r, 0) \text{ and } \exists \lim_{s \to 0^-} \psi(s) \in \mathbf{R}^n\}$$

is considered in section 8.2 of [4]. Its elements are $\psi = \varphi + X_0 \sigma$, with $\varphi \in B$, $\sigma \in \mathbf{R}^n$ (column vector) and

$$X_0(s) = \begin{cases} 0, & -r \leq s < 0 \\ I_n, & s = 0, \end{cases}$$

where $0, I_n$ are the zero and, respectively, the unity nxn matrix. The norm of $\psi$ is defined as the sum of the norm of $\varphi$ in B and the norm of $\sigma$ in $\mathbf{R}^n$. The complexifications $B_\mathbf{C}, B_{0\mathbf{C}}$ of B, respectively $B_0$, are used below.

Consider $\delta(\cdot)$ - the Dirac function, and the nxn matrix valued function $\Delta(s) = \delta(s) I_n$. Also consider the bounded linear operator $L_\alpha : B \to \mathbf{R}^n$,

$$L_\alpha \varphi = \int_{-r}^{0} d\eta_\alpha(s) \varphi(s),$$

with $\eta_\alpha(s) = A(\alpha)\Delta(s) - B(\alpha)\Delta(s+r)$. By denoting $x_t(s) = x(t+s)$, $s \in [-r,0]$, we have in the spirit of [4], the following relations, equivalent with (1), (2):

$$\frac{dx_t}{dt}(0) = L_\alpha(x_t) + f(x_t(0), x_t(-r), \alpha), \tag{3}$$

$$\frac{dx_t}{dt}(s) = \frac{dx_t}{ds}(s), \tag{4}$$

$$x_0 = \varphi. \tag{5}$$

Define (see also [4]) the linear operator, $\tilde{A}_\alpha \varphi = \dot\varphi + X_0\left[L_\alpha(\varphi) - \dot\varphi(0)\right]$,

$\tilde{A}_\alpha : C^1([-r,0], \mathbf{R}^n) \subset BC \to BC$. Now we can rewrite the above problem as

$$\frac{dx_t}{dt} = \tilde{A}_\alpha x_t + X_0 f(x_t(0), x_t(-r), \alpha), \tag{6}$$

$$x_0 = \varphi. \tag{7}$$

The last term of (6) may be written as $X_0 f(\int_{-r}^0 d\Delta(s)x_t(s), \int_{-r}^0 d\Delta(s+r)x_t(s), \alpha)$. We define $F(x_t, \alpha) = f(\int_{-r}^0 d\Delta(s)x_t(s), \int_{-r}^0 d\Delta(s+r)x_t(s), \alpha)$. Thus (6) and (7) take the form

$$\frac{dx_t}{dt} = \tilde{A}_\alpha x_t + X_0 F(x_t, \alpha), \tag{8}$$

$$x_0 = \varphi, \tag{9}$$

this being the abstract problem in $B_0$ equivalent to (1),(2).

The eigenvalues of its infinitesimal generator are the roots of the equation
$$\det(\lambda I - A(\alpha) - e^{-\lambda r} B(\alpha)) = 0.$$
We assume the following hypothesis, that we denote *H1*.

*H1. An open set U exists in the parameter plane such that for every $\alpha \in U$, there is a pair of complex conjugated simple eigenvalues $\lambda_{1,2}(\alpha) = \mu(\alpha) \pm i\omega(\alpha)$, with the property that there is a $\alpha_0 \in U$ such that $\lambda_{1,2}(\alpha_0) = \pm i\omega(\alpha_0) = \pm i\omega_0$, with $\omega_0 > 0$ and for every $\alpha \in U$, all other eigenvalues have strictly negative real parts, uniformly bounded from above by a negative number.*

By a simple eigenvalue we mean an eigenvalue having the algebraic multiplicity equal to 1.

We remark that *H1* implies the existence of a neighborhood of $\alpha_0$ such that each eigenvalue different from $\lambda_{1,2}(\alpha)$ has real part strictly less than $\mu(\alpha)$.

The eigenfunctions corresponding to $\lambda_i(\alpha)$, $i=1,2$, are elements of $B_\mathbf{C}$ -the complexification of $B$, namely

$$\varphi_i(\alpha)(s) = e^{\lambda_i(\alpha)s} \varphi_i(\alpha)(0), \quad s \in [-r,0],$$

where $\varphi_i(\alpha)(0)$ is a solution of $(\lambda_i(\alpha)I - A(\alpha) - e^{-\lambda_i(\alpha)r} B(\alpha))\varphi = 0$. Obviously, $\varphi_2(\alpha) = \overline{\varphi_1(\alpha)}$.

Denote by $\mathbf{M}_{\{\lambda_{1,2}(\alpha)\}}$ the linear subspace of $B_\mathbf{C}$, spanned by $\{\varphi_1(\alpha), \varphi_2(\alpha)\}$ and $\Phi(\alpha)$ the matrix having as columns the vectors $\varphi_1(\alpha), \varphi_2(\alpha)$. Let $\{\psi_1(\alpha), \psi_2(\alpha)\}$ be two eigenvectors for the adjoint problem ([3], [4]), corresponding to the eigenvalues $-\lambda_{1,2}(\alpha)$ of the infinitesimal generator of the adjoint problem. They are elements of $B_\mathbf{C}^*$ - the complexification of $B^* = C([0,r], \mathbf{R}^n)$ -row vectors, and we assume that they are selected such that, $\Psi(\alpha)$ being the matrix having as rows the vectors $\psi_1(\alpha), \psi_2(\alpha)$, the relation $(\Psi(\alpha), \Phi(\alpha)) = I_2$ holds, where $(\cdot,\cdot): B_C^* \times B_C \to C$ is defined by

$$(\psi, \varphi) = \psi(0)\varphi(0) - \int_{-r}^{0}\int_{0}^{\theta} \psi(\xi - \theta) d\eta_\alpha(\theta)\varphi(\xi)d\xi, \quad \psi \in B_C^*, \varphi \in B_C.$$

In [4] a projection $\pi: B_{0C} \to \mathbf{M}_{\{\lambda_{1,2}(\alpha)\}}$, is defined by $\pi(\varphi + X_0 a) = \Phi(\alpha)[(\Psi(\alpha), \varphi) + \Psi(0)a]$.
With this projection the space $B_{0C}$ is decomposed as $B_{0C} = \mathbf{M}_{\{\lambda_{1,2}(\alpha)\}} \oplus \operatorname{Ker}\pi$. Since the solution $x_t$ of (8) and (9) belongs to $C^1([-r, 0], \mathbf{R}^n)$, it is decomposed as

$$x_t = \Phi(\alpha)u(t) + y(t), \tag{10}$$

with $u(t) = (\Psi(\alpha), x_t)$ and $y(t) = (I - \pi)x_t$, where $u(t) = (z(t), \bar{z}(t))$ -column vector, $z(t) \in C$.
Let us define

$$B(\alpha) = \begin{pmatrix} \lambda_1(\alpha) & 0 \\ 0 & \bar{\lambda}_1(\alpha) \end{pmatrix}.$$

The projection of eq. (8) on $\mathbf{M}_{\{\lambda_{1,2}(\alpha)\}}$ is

$$\Phi(\alpha)\dot{u} = \Phi(\alpha)B(\alpha)u + \Phi(\alpha)\Psi(\alpha)(0)F(\Phi(\alpha)u(t) + y(t), \alpha),$$

and since $\Phi(\alpha)$ is invertible, this is equivalent to

$$\dot{u} = B(\alpha)u + \Psi(\alpha)(0)F(\Phi(\alpha)u(t) + y(t), \alpha). \tag{11}$$

By projecting the initial condition we find

$$u(0) = (\Psi(\alpha), \varphi).$$

**3. Existence of the invariant manifold and the restricted equation**

If $\alpha \in U \setminus \{\alpha_0\}$ and, $\operatorname{Re} \lambda_{1,2}(\alpha) > 0$ then, since these are the only two eigenvalues with positive real part (and they are simple), there is a local invariant manifold for the problem, namely the local unstable manifold, tangent to the space $\mathbf{M}_{\{\lambda_{1,2}(\alpha)\}}$, [3], [7].

For $\alpha = \alpha_0$, since $\operatorname{Re} \lambda_{1,2}(\alpha_0) = 0$, there is a local invariant manifold for the problem, namely the local center manifold, tangent to the space $\mathbf{M}_{\{\lambda_{1,2}(\alpha_0)\}}$, [7].

Hence, for every $\alpha \in U$ with $\mu(\alpha) \geq 0$, there is a neighborhood $V(\alpha)$ of $0 \in B$, and a local invariant manifold $W_{loc}(\alpha) \subset V(\alpha)$, which is the graph of a $C^1$ function. That is, the local invariant manifold may be expressed as

$$W_{loc}(\alpha) = \{\varphi + w_\alpha(\varphi); \varphi \in \mathbf{M}_{\{\lambda_{1,2}(\alpha)\}} \cap V(\alpha)\},$$

where $w_\alpha: \mathbf{M}_{\{\lambda_{1,2}(\alpha)\}} \to \operatorname{Ker}\pi$ is a $C^1$ function, $w_\alpha(0) = 0$, and it has zero differential at 0. Since $\varphi \in \mathbf{M}_{\{\lambda_{1,2}(\alpha)\}}$, we have $\varphi(\alpha) = z\varphi_1(\alpha) + \bar{z}\bar{\varphi}_1(\alpha)$, with $z \in C$. This relation induces a dependence of $z, \bar{z}$ to $w_\alpha(\varphi)$ that justifies the notation $w_\alpha(\varphi) = w_\alpha(z, \bar{z})$.
Equation (11) implies

$$\dot{z}(t) = \lambda_1(\alpha)z(t) + \psi_1(\alpha)(0)F(z(t)\varphi_1(\alpha) + \bar{z}(t)\bar{\varphi}_1(\alpha) + w_\alpha(z(t), \bar{z}(t)), \alpha), \tag{12}$$

$$z(0) = (\psi_1(\alpha), \varphi). \tag{13}$$

Let $S_\alpha(t)\phi$ be the solution of eq. (1) corresponding to the initial condition $\phi$, at the moment $t$. If $\phi \in W_{loc}(\alpha)$, then

$$S_\alpha(t)\phi = z(t)\varphi_1(\alpha) + \bar{z}(t)\bar{\varphi}_1(\alpha) + w_\alpha(z(t),\bar{z}(t)).\tag{14}$$

By using again the function $f$, (12) becomes

$$\dot{z}(t) = \lambda_1(\alpha)z(t) + \psi_1(\alpha)(0)f(S_\alpha(t)\phi(0), S_\alpha(t)\phi(-r), \alpha),\tag{15}$$

or,

$$\dot{z}(t) = \lambda_1(\alpha)z(t) + g(z(t),\bar{z}(t),\alpha),\tag{16}$$

by denoting

$$g(z(t),\bar{z}(t),\alpha) = \psi_1(\alpha)(0)f(S_\alpha(t)\phi(0), S_\alpha(t)\phi(-r), \alpha).\tag{17}$$

## 4. The equations for the invariant manifold

The following proposition is a natural consequence of the invariance of $W_{loc}(\alpha)$. A similar result is given in [7], on the center manifold. We give the proof for the sake of completeness.

**Proposition 1.** *Let $\phi \in W_{loc}(\alpha)$ be the initial value for the problem (1). Then the function $w_\alpha$ satisfies the following equations*

$$\frac{\partial}{\partial t}w_\alpha(z(t),\bar{z}(t))(s) + g(z(t),\bar{z}(t),\alpha)\varphi_1(\alpha)(s) + \bar{g}(z(t),\bar{z}(t),\alpha)\bar{\varphi}_1(\alpha)(s) =$$
$$= \frac{\partial}{\partial s}w_\alpha(z(t),\bar{z}(t))(s), \qquad s \in [-r,0],\tag{18}$$

$$\frac{\partial}{\partial t}w_\alpha(z(t),\bar{z}(t))(0) + g(z(t),\bar{z}(t),\alpha)\varphi_1(\alpha)(0) + \bar{g}(z(t),\bar{z}(t),\alpha)\varphi_2(\alpha)(0) =$$
$$= A(\alpha)w_\alpha(z(t),\bar{z}(t))(0) + B(\alpha)w_\alpha(z(t),\bar{z}(t))(-r) + f(S_\alpha(t)\phi(0), S_\alpha(t)\phi(-r), \alpha)\tag{19}$$

*with $z(t)$ solution of the Cauchy problem (16),(13) and $g$ defined by (17).*

**Proof.** Since $\phi \in W_{loc}(\alpha)$ and $W_{loc}(\alpha)$ is invariant, $S_\alpha(t)\phi \in W_{loc}(\alpha)$.
Let us denote, for $t \geq 0$ and $s \in [-r,0]$, $S_\alpha(t)\phi(s) = x(\phi)(t+s)$.
Obviously

$$\frac{\partial x(\phi)}{\partial t}(t+s) = \frac{\partial x(\phi)}{\partial s}(t+s).$$

This and (14) imply

$$\frac{\partial}{\partial t}w_\alpha(z(t),\bar{z}(t))(s) + \dot{z}(t)\varphi_1(\alpha)(s) + \dot{\bar{z}}(t)\bar{\varphi}_1(\alpha)(s) =$$
$$= \frac{\partial}{\partial s}w_\alpha(z(t),\bar{z}(t))(s) + z(t)\dot{\varphi}_1(\alpha)(s) + \bar{z}(t)\dot{\bar{\varphi}}_1(\alpha)(s)\tag{20}$$

(here $\dot{\varphi}_1(\alpha)(s) = \frac{d}{ds}\varphi_1(\alpha)(s)$) and thus

$$\frac{\partial}{\partial t}w_\alpha(z(t),\bar{z}(t))(s) + \left[\dot{z}(t) - \lambda_1(\alpha)z(t)\right]\varphi_1(\alpha)(s) + \left[\dot{\bar{z}}(t) - \overline{\lambda_1(\alpha)}\bar{z}(t)\right]\bar{\varphi}_1(\alpha)(s) = \frac{\partial}{\partial s}w_\alpha(z(t),\bar{z}(t))(s).\tag{21}$$

With (16) we obtain (18).
On another side, since $S_\alpha(t)\phi$ is a solution of equation (1), we have

$$\dot{z}(t)\varphi_1(\alpha)(s) + \dot{\bar{z}}(t)\overline{\varphi_1}(\alpha)(s) + \frac{\partial}{\partial t}w_\alpha(z(t),\bar{z}(t))(s) = A(\alpha)S_\alpha(t)\phi(s) + B(\alpha)S_\alpha(t-r)\phi(s) +$$
$$+ f(S_\alpha(t)\phi(s), S_\alpha(t-r)\phi(s)),$$
and, by taking $s = 0$, we obtain (19).

This proposition allows the determination of the coefficients of the series of powers in $z$ and $\bar{z}$ of the function $w_\alpha$. Indeed, let us write

$$f(S_\alpha(t)\phi(0), S_\alpha(t)\phi(-r),\alpha) = \sum_{j+k\geq 2} \frac{1}{j!k!} F_{jk}(\alpha) z^j \bar{z}^k, \tag{22}$$

$$g(z(t),\bar{z}(t),\alpha) = \sum_{j+k\geq 2} \frac{1}{j!k!} g_{jk}(\alpha) z^j \bar{z}^k, \quad w_\alpha(\theta,z,\bar{z}) = \sum_{j+k\geq 2} \frac{1}{j!k!} w_{jk}(\theta,\alpha) z^j \bar{z}^k \tag{23}$$

where $g_{jk}(\alpha) = \psi_1(0) F_{jk}(\alpha)$.

By replacing (22) and (23) in (18) and by matching the obtained series, we get first order linear differential equations for $w_{jk}$.

Thus, equation (18) implies

$$\sum_{j+k\geq 2} \frac{1}{j!k!} \frac{d}{ds} w_{jk}(s,\alpha) z^j \bar{z}^k = \left( \sum_{j+k\geq 2} \frac{1}{j!k!} g_{jk}(\alpha) z^j \bar{z}^k \right) \varphi_1(\alpha)(s) +$$
$$+ \left( \sum_{j+k\geq 2} \frac{1}{j!k!} \bar{g}_{jk}(\alpha) \bar{z}^j z^k \right) \overline{\varphi_1}(\alpha)(s) + \sum_{j+k\geq 2} \frac{1}{j!k!} w_{jk}(s,\alpha) \left( jz^{j-1}\bar{z}^k \dot{z} + kz^j \bar{z}^{k-1} \dot{\bar{z}} \right). \tag{24}$$

In this equality $\dot{z}$ and $\dot{\bar{z}}$ will be replaced with the right hand side of (16) to obtain

$$\sum_{j+k\geq 2} \frac{1}{j!k!} \frac{d}{ds} w_{jk}(s,\alpha) z^j \bar{z}^k = \left( \sum_{j+k\geq 2} \frac{1}{j!k!} g_{jk}(\alpha) z^j \bar{z}^k \right) \varphi_1(\alpha)(s) +$$
$$+ \left( \sum_{j+k\geq 2} \frac{1}{j!k!} \bar{g}_{jk}(\alpha) \bar{z}^j z^k \right) \overline{\varphi_1}(\alpha)(s) + \sum_{j+k\geq 2} \frac{1}{j!k!} w_{jk}(s,\alpha) (j\lambda_1(\alpha) + k\overline{\lambda_1}(\alpha)) z^j \bar{z}^k + \tag{25}$$
$$+ \sum_{j+k\geq 2} \frac{1}{j!k!} w_{jk}(s,\alpha) \left( jz^{j-1}\bar{z}^k \left( \sum_{l+m\geq 2} \frac{1}{l!m!} g_{lm}(\alpha) z^l \bar{z}^m \right) + kz^j \bar{z}^{k-1} \left( \sum_{l+m\geq 2} \frac{1}{l!m!} \bar{g}_{lm}(\alpha) \bar{z}^l z^m \right) \right).$$

By matching the same order terms we obtain first order differential equations for $w_{jk}(.,\alpha)$.

A relation similar to (25) is obtained by substituting the series (22), (23) in (19), and by using (16):

$$\sum_{j+k\geq 2} \frac{1}{j!k!} w_{jk}(0,\alpha)(j\lambda_1(\alpha) + k\overline{\lambda_1}(\alpha)) z^j \bar{z}^k +$$
$$+ \sum_{j+k\geq 2} \frac{1}{j!k!} w_{jk}(0,\alpha) \left( jz^{j-1}\bar{z}^k \left( \sum_{l+m\geq 2} \frac{1}{l!m!} g_{lm}(\alpha) z^l \bar{z}^m \right) + kz^j \bar{z}^{k-1} \left( \sum_{l+m\geq 2} \frac{1}{l!m!} \bar{g}_{lm}(\alpha) \bar{z}^l z^m \right) \right) +$$
$$+ \left( \sum_{j+k\geq 2} \frac{1}{j!k!} g_{jk}(\alpha) z^j \bar{z}^k \right) \varphi_1(\alpha)(0) + \left( \sum_{j+k\geq 2} \frac{1}{j!k!} \bar{g}_{jk}(\alpha) \bar{z}^j z^k \right) \overline{\varphi_1}(\alpha)(0) = \tag{26}$$
$$= A(\alpha) \sum_{j+k\geq 2} \frac{1}{j!k!} w_{jk}(0,\alpha) z^j \bar{z}^k + B(\alpha) \sum_{j+k\geq 2} \frac{1}{j!k!} w_{jk}(-r,\alpha) z^j \bar{z}^k + \sum_{j+k\geq 2} \frac{1}{j!k!} F_{jk}(\alpha) z^j \bar{z}^k.$$

The relations obtained by equating the terms with similar powers of $z, \bar{z}$ in (26) are used as conditions to determine the constants that appear in the general form of $w_{jk}$ obtained above.

We firstly remark that the coefficients of the second order terms in $z$ and $\bar{z}$ in the expansion of $f_i(S_\alpha(t)\phi(0), S_\alpha(t)\phi(-r), \alpha)$, $i = 1,...,n$, are independent on those of $w_\alpha$, they depend only on the coefficients of the Taylor series of $f_i(x, y, \alpha)$. The similar assertion holds for the coefficients of $g(z(t), \bar{z}(t), \alpha)$. Hence $g_{20}(\alpha), g_{11}(\alpha), g_{02}(\alpha)$ are known, given $f_i$ and $\psi_1(\alpha)$.

The following algorithm to determine $w_{jk}(\alpha)$ must be used.
- $w_{20}(\alpha), w_{11}(\alpha), w_{02}(\alpha)$ are determined from the equations obtained by identifying the terms containing $z^2, z\bar{z}, \bar{z}^2$ respectively, in (25), with initial conditions obtained by the same method from (26); they depend on $g_{20}(\alpha), g_{11}(\alpha), g_{02}(\alpha)$.
- $w_{20}(\alpha), w_{11}(\alpha), w_{02}(\alpha)$ are used to compute $g_{jk}(\alpha), j+k=3$, (from (14), (17), (23) ).
- $w_{jk}(\alpha), j+k=3$, are determined from the equations (25) and conditions (26); they depend on $g_{jk}(\alpha), j+k \leq 3$.
- $w_{jk}(\alpha), j+k \leq 3$ are used to compute $g_{jk}(\alpha), j+k=4$ (from (14), (17), (23) ).
- $w_{jk}(\alpha), j+k=4$, are determined from the equations (25) and conditions (26); they depend on $g_{jk}(\alpha), j+k \leq 4$.
- $w_{jk}(\alpha), j+k \leq 4$ are used to compute $g_{jk}(\alpha), j+k=5$.

We do not need so many terms to have an accurate form of the invariant manifold, but we need them in order to discuss the behaviour of the solution $z$ of (16) that determines the solution of (1) on the invariant manifold (see below).

## 5. The Bautin bifurcation

After applying the above algorithm, equation (16) has the form

$$\dot{z}(t) = \lambda_1(\alpha)z(t) + \sum_{2 \leq j+k \leq 5} \frac{1}{j!k!} g_{jk}(\alpha) z^j \bar{z}^k + O(|z|^6). \tag{27}$$

Now, let $l_1(\alpha), l_2(\alpha)$ be the first and the second Lyapunov coefficient respectively, associated to (27), defined in [6]. They are functions of $g_{kl}(\alpha)$. More specific, $l_1(\alpha)$ is function of $g_{kl}(\alpha)$ with $2 \leq k+l \leq 3$, and $l_2(\alpha)$ is function of $g_{kl}(\alpha)$ with $2 \leq k+l \leq 5$.

Let us define, with [6], $v_1 = \dfrac{\mu(\alpha)}{\omega(\alpha)}$, $v_2 = l_1(\alpha)$, $v = (v_1, v_2)$.

Let us consider the following hypothesis, denoted *H2*.

**H2.** $l_1(\alpha_0) = 0$, $l_2(\alpha_0) \neq 0$, and the map $(\alpha_1, \alpha_2) \to (v_1, v_2)$ is regular at $\alpha_0$.

The value in $\alpha_0$ of the first Lyapunov coefficient is [6]

$$l_1(\alpha_0) = \frac{1}{2\omega_0^2} \text{Re}(ig_{20}(\alpha_0)g_{11}(\alpha_0) + \omega_0 g_{21}(\alpha_0)).$$

The value in $\alpha_0$ of the second Lyapunov coefficient is much more complicated and we do not reproduce it here (see [6]). We only note that, since the expression of $l_2(\alpha_0)$ contains coefficients

$g_{ij}$ with $i + j \le 5$, in order to compute it, we have to determine the functions $w_{ij}(.)$ with $i + j \le 4$, as the algorithm presented above shows.

Th. 8.2. of [6] asserts that, if *H2* is satisfied for eq. (16), then eq. (16) may be transformed, by smooth invertible changes of the complex function *z* and time reparametrizations, into

$$\dot{z} = (v_1 + i)z + v_2 z|z|^2 + L_2(v)z|z|^4 + O(|z|^6), \tag{28}$$

where $L_2(v) = l_2(\alpha(v))$.
We assume that *H2* is satisfied for our problem.
Problem (28) is then (see [6]) locally, near *z*=0, topologically equivalent to the following complex normal form of the Bautin bifurcation

$$\dot{z} = (\beta_1 + i)z + \beta_2 z|z|^2 + s\, z|z|^4, \tag{29}$$

where $\beta_1 = v_1, \beta_2 = \sqrt{|L_2(v)|}v_2$, $s = \text{sign}\, l_2(\alpha_0)$.

Assume the third important fact denoted by *H3*.
*H3.* $l_2(\alpha_0) > 0$ (hence *s*=1).

Consider the polar form of (29), in the hypothesis *H3*

$$\begin{cases} \dot{\rho} = \rho(\beta_1 + \beta_2 \rho^2 + \rho^4), \\ \dot{\eta} = 1, \end{cases} \tag{30}$$

with $(\rho, \eta)$ the polar co-ordinates.
Due to *H2*, there is a $C^1$ bijection *T* between a neighbourhood $U_1$ of $\alpha_0$, in the $(\alpha_1, \alpha_2)$ plane and a neighbourhood *W* of (0,0) in the $(\beta_1, \beta_2)$ plane.

We study the behaviour of the solutions only for $\beta_1 \ge 0$ ( $\beta_1$ and $\mu(\alpha)$ have the same sign and we proved the existence of the invariant manifold only for $\mu(\alpha) \ge 0$ ).

For those $\beta \in W$ with $\beta_1 \ge 0, \beta_2 \ge 0$, (30) has an repulsive focus in 0, and so does (29). Then, since the dynamical system generated by (29) is topologically equivalent to that generated by (16), this one has an unstable equilibrium in 0 (a repulsive focus or node, since they cannot be distinguished by topological equivalence).
In order to transport these conclusions to the solution of eq. (1), we remind relation (14):

$$S_\alpha(t)\phi(\theta) = z(t)\varphi_1(\alpha)(\theta) + \bar{z}(t)\overline{\varphi_1}(\alpha)(\theta) + w_\alpha(\theta, z(t), \bar{z}(t)).$$

It follows that for the corresponding (by $T^{-1}$) zone of $\alpha$ plane, let us denote it $V_1$, equation (1) has 0 as an unstable equilibrium point (focus or node) on the invariant manifold. Due to the relation between $\beta_1$ and $\mu(\alpha)$, this holds on the unstable manifold for $\mu(\alpha) > 0$ and, for $\alpha_0$ (where $\mu(\alpha_0) = 0$) on the center manifold.

In the intersection of the quadrant $\beta_1 \ge 0, \beta_2 < 0$ and *W* there is a zone, situated between the axis $\beta_2 = 0$ and the curve $\Gamma = \{(\beta_1, \beta_2); \beta_2 = -2\sqrt{\beta_1}\}$ where eq. $\beta_1 + \beta_2 \rho^2 + \rho^4 = 0$ has no positive solutions and thus there are no limit cycles for (30). The local phase portrait for (1) on the invariant manifold is as described above (unstable equilibrium point), for the corresponding (by $T^{-1}$) zone in the $\alpha$ plane, let us denote it by $V_2$. Now let us put $V = V_1 \cup V_2$.

In the zone situated between the curve $\Gamma = \{(\beta_1, \beta_2); \beta_2 = -2\sqrt{\beta_1}\}$ and the axis $\beta_1 = 0$ eq. $\beta_1 + \beta_2 \rho^2 + \rho^4 = 0$ has two positive solutions, $\rho_1(\beta), \rho_2(\beta)$. Thus two concentric closed orbits for

(30), hence for (29), exist. Since, for the corresponding zone of the $\alpha$ plane (let us denote it by $V^*$), the problems (29) and (16) are topologically equivalent, eq. (16) will also have two periodic solutions. Relation (14) shows that, in this case, on the unstable manifold there exist two limit cycles for (1).

As we cross $\Gamma$, leaving the zone in the $\beta$ plane described above, the two limit cycles collide and disappear.

## 6. Conclusions

The facts discussed above lead to the following result.

**Theorem** *If H1, H2, H3 are satisfied for eq. (1), then at $\alpha_0$ a Bautin type bifurcation takes place. There is a neighbourhood $U_1$ of $\alpha_0$ in the $\alpha$ plane having a subset $V$ (with $\alpha_0 \in V$) with the property that for every $\alpha \in V$, 0 is an unstable equilibrium point (focus or node) for the problem (1) restricted to the bidimensional invariant manifold defined above.*

*There is also a subset $V^*$ of $U_1$, (having $\alpha_0$ as a limit point) with the property that for every $\alpha \in V^*$, the restriction of problem (1) to the unstable manifold has two limit cycles (one inside the other). In this case also* **0** *is an unstable equilibrium point.*

Let us remark that the interior limit cycle is attractive, while the exterior limit cycle is repulsive.